Article: CJB/2010/3

# Generalizations of Hagge's Theorems

**Christopher J Bradley**

**Abstract:** Two generalizations of Hagge's theorems are described. In the first we consider what happens when one moves from the orthocentre to a general point. What one loses by doing so is the indirect similarity and hence one loses the centre of indirect similarity. Instead one proceeds from the centre of the circle under consideration. In the second generalization we consider pairs of triangles that have orthologic centres with respect to each other, so that an indirect similarity is the main feature preserved.

## 1. Introduction

There are two unrelated generalizations of the Hagge configuration.

The first generalization involves moving from the orthocentre *H* to another point *K* at the expense of losing the indirect similarity. It also appears that the material of Article: CJB/2010/2 cannot be completely generalized, as there is a restriction on the facility for creating four circles from a cyclic quadrilateral. Also it does not appear there is a simple prescription that relates the point of perspective *P* between the triangles *XYZ* and *UVW* and the centre *Q* of the circle through *K* on which they lie. In the case of the Hagge construction Peiser [1] showed that the isogonal conjugate *Pg* of *P* is located in such a way that the nine-point centre *N* is the midpoint of *QPg*. In Sections 3 – 6 we give proofs of the main results that form the extensions of Hagge's theorems and the four Hagge circle property. In the Hagge construction *AP*, *BP*, *CP* are drawn to meet the circumcircle *Γ* at points *D*, *E*, *F* respectively. These points are then reflected in the sides *BC*, *CA*, *AB* to create three new points *U*, *V*, *W*. The first Hagge theorem is that the circle *UVW* passes through the orthocentre *H*. Next one draws *AH*, *BH*, *CH* to meet this circle at *X*, *Y*, *Z* respectively and the second Hagge theorem is that *UX*, *VY*, *WZ* are concurrent at *P*. It is interesting that in the present generalization one does not start with a point *P*, but with the centre *Q* of the proposed circle *UVWK*.

The second generalization is designed to maintain the property that triangles *ABC* and *XYZ* have orthologic points with respect to each other and this leads to these triangles and their circumcircles being related by an indirect similarity and also to triangles *DEF* and *UVW* being related by the same indirect similarity. What is lost is the reflection property mentioned in the last paragraph. We describe these generalizations in turn and the following results hold for the first:

## 2. Results for the first generalization

**Theorem 1**

Let $ABC$ be a triangle and $K$ any point not lying on its sides or extensions. With centre any point $Q$ draw a circle $\Sigma$ to pass through $K$. The intersections with $\Sigma$ of circles $BKC$, $CKA$, $AKB$ are denoted by $U$, $V$, $W$ respectively. Then

$$(BU/CU)(CV/AV)(AW/BW) = -1. \tag{2.1}$$

In this expression $(BU/CU)$, for example, is taken as positive if $KBUC$ is convex and negative otherwise. The converse is also true, that if $U$, $V$, $W$ lie on a circle and this relation holds, then the circle passes through $K$.

This result is due to Boreico [2].

**Theorem 2**

Given the configuration of Theorem 1, suppose now $AK$, $BK$, $CK$ meet $\Sigma$ at $X$, $Y$, $Z$ respectively, then $UX$, $VY$, $WZ$ are concurrent at a point $P$.

In other words the pairs $(U, X)$, $(V, Y)$, $(W, Z)$ are in involution on $\Sigma$, by means of a perspective point. These results are shown in Fig. 1.

## 3. Proof of Theorem 1

Invert with respect to $K$ and denote inverse points by primes. We require $U'$, $V'$, $W'$ collinear, which is the case by Menelaus' Theorem if, and only if, $(B'U'/C'U')(C'V'/A'V')(A'W'/B'W') = -1$. However, $B'U'/C'U' = (BU/CU)(KB'/KC')$. Multiplying three such relationships we get the required result.

## 4. Proof of Theorems 2, 3 and 4

**Theorem 2**

We use Cartesian co-ordinates with $K$ as origin and let $\Sigma$ have centre $Q$ with co-ordinates $(-\tfrac{1}{2}, -\tfrac{1}{2})$. This merely imposes a scale and a direction on the configuration, so there is no loss of generality. Then $\Sigma$

$$x^2 + y^2 + x + y = 0. \tag{4.1}$$

We parameterize $\Sigma$ by taking lines through $K$ to be of the form $y = mx$ and then $m$ serves as the parameter of the point where this line meets $\Sigma$ again.

Circle $BKC$ and $\Sigma$ meet on the line $KU$ with equation $(p - 1)x + (q - 1)y = 0$, so the parameter for the point $U$ is $a = -(p - 1)/(q - 1)$. Similarly the parameters for $V$ and $W$ are respectively $c = -(s - 1)/(t - 1)$ and $e = -(u - 1)/(v - 1)$. Now circles $CKA$ and $AKB$ meet at $A$ and $K$, so their common chord is the line $AK$, which therefore has equation $(s - u)x + (t - v)y = 0$. It follows that the parameter of $X$ is $b = -(s - u)/(t - v)$. Similarly $Y$ and $Z$ have parameters that are respectively $d = -(p - u)/(q - v)$ and $e = -(p - s)/(q - t)$.

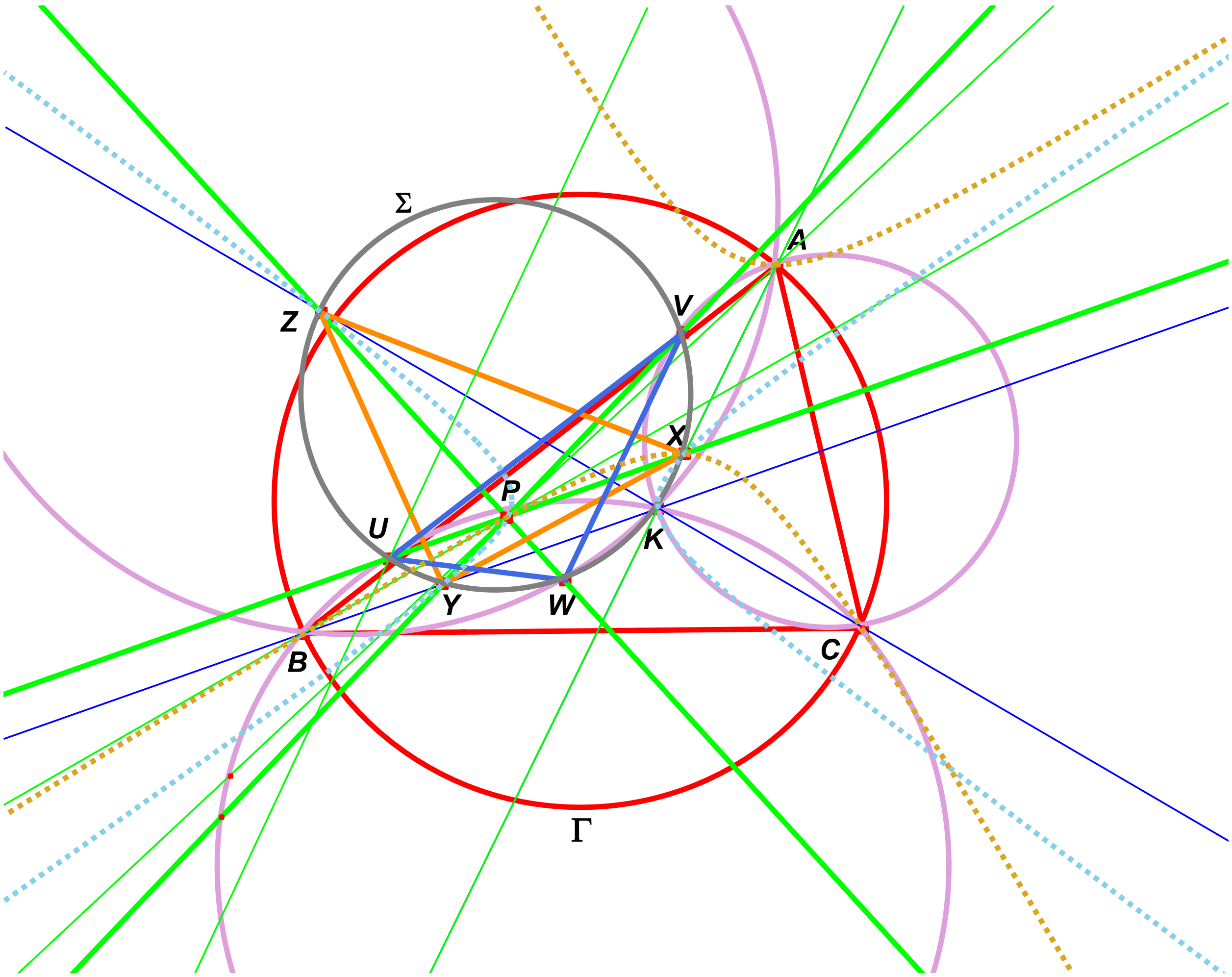

**Fig. 1**
**The involution on Σ**

To prove that $UX$, $VY$, $WZ$ are concurrent it is sufficient to show that the pairs $(a, b)$, $(c, d)$, $(e, f)$ are in involution. This is because pairs in involution on a conic must arise from a vertex of perspective. Since an involution of pairs $(h, k)$ is specified by an equation of the form $lhk + m(h + k) + n = 0$ for some real numbers $l$, $m$, $n$ ($m^2 \neq nl$) it follows that the pairs are in involution if the determinant with rows$(ab, a + b, 1)(cd , c + d, 1)$, $(ef, e + f, 1)$ vanishes. *DERIVE* verifies this is indeed the case.

**Theorem 3**
Let *VW* meet *AKX* at *L*, with *M*, *N* similarly defined, then *L*, *M*, *N*, *P* are collinear.

**Proof**
Consider the hexagon *VWUXKY* on the circle *Σ*. *VW*^*XK* = *L*, *WU*^*KY* = *M* and *UX*^*YV* = *P*. It follows by Pascal's theorem that *L*, *M*, *P* are collinear. Similarly *M*, *N*, *P* are collinear.

The above theorems may be generalized even further, as the following theorem shows:

**Theorem 4**
Let *ABC* be a triangle and *D*, *E*, *F* any three generally situated points. Draw any conic *Σ* through *D*, *E*, *F*. Let conic *BCDEF* meet *Σ* at *U*, conic *CADEF* meet *Σ* at *V* and conic *ABDEF* meet *Σ* at *W*. Further let *AD*, *BD*, *CD* meet *Σ* at *X*, *Y*, *Z* respectively, then *UX*, *VY*, *WZ* are concurrent.

**Proof**
Project *E*, *F* to the circular points at infinity.
Theorem 4 is illustrated in Fig. 2

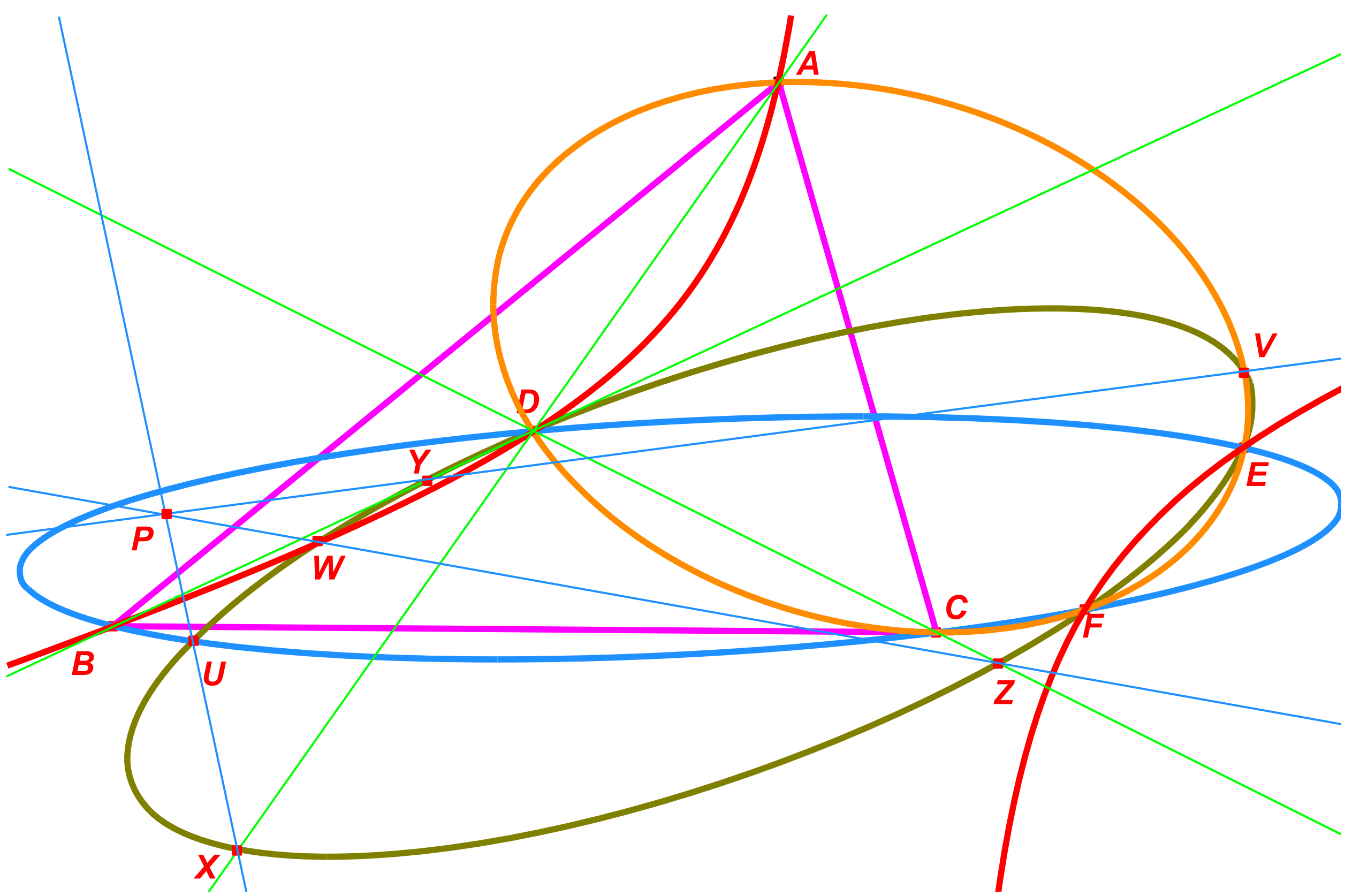

**Fig.2**
**Illustration of Theorem 4**

## 5. Restricted generalization of the four Hagge circle configuration

We refer back to Article CJB/2010/2 for a detailed description of this configuration, but in brief what happens is as follows: Take a triangle $ABC$ inscribed in a circle $\Gamma$. Let its orthocentre be $D_4$. Take a point $P$ anywhere (not on $\Gamma$ or the sides of $ABC$), that will serve as a centre of inverse similarity. Draw the rectangular hyperbola through $A$, $B$, $C$, $D_4$ and $P$. Its centre, as is well known, is a point $M$ on the nine-point circle of $ABC$. Now perform a 180° rotation of $A$, $B$, $C$, $D_4$ about $M$ to get the points $A_1$, $B_2$, $C_3$, $D$. Then, wherever the initial selection of $P$ is made, it is always the case that $D$ lies on $\Gamma$ and the points $A_1$, $B_2$, $C_3$ are the orthocentres of triangles $BCD$, $ACD$, $ABD$ respectively. The four Hagge circles of $P$ with respect to the triangles $BCD$, $ACD$, $ABD$, $ABC$ can now be drawn using the same inverse spiral symmetry about $P$ in each case.

If we now choose $D_4$ to be a general point of triangle $ABC$ other than the orthocentre, the same construction does not work. If $P$ is chosen anywhere (not on $\Gamma$ or the sides of $ABC$), then it turns out that the point $D$ does not lie on $\Gamma$. In fact there is only one conic through $A$, $B$, $C$, $D_4$ for which this turns out to be true. The point $M$ that is the centre of this conic lies somewhere on the circular locus of points for which $D_4M$ produced to meet $\Gamma$ at $D$ is such that $D_4M = MD$. $P$ is restricted to lie anywhere on this conic, which may be a hyperbola or an ellipse.

The next question to ask is how to draw the four circles that have the point $P$ as their involution point. As we no longer have either Peiser's [1] prescription or the inverse spiral symmetry to rely upon, there has to be some other method of obtaining the four circle centres from $P$. We draw on the property of the four Hagge circle configuration to provide the answer. In that configuration it turned out that points $A_k$ ($k$ = 1, 2, 3, 4) lie on a line through $P$. We call this the $A$-line. Similarly there is a $B$-line, a $C$-line and a $D$-line. We insist for the generalization that these lines must exist. We now describe how to obtain the points $A_4$, $B_4$, $C_4$. $A_4$ is the intersection of $A_1P$ with $AD_4$, $B_4$ is the intersection of $B_2P$ with $BD_4$ and $C_4$ is the intersection of $C_3P$ with $CD_4$. The circle $\Sigma_4$ is now defined to be the circle $A_4B_4C_4$.

In order to justify the validity of the construction certain theorems have to be proved.

**Theorem 5**
The circle $\Sigma_4$ passes through $D_4$ and similarly $\Sigma_1$, $\Sigma_2$, $\Sigma_3$ pass through $A_1$, $B_2$, $C_3$

**Theorem 6**
Let circle $BCD_4$ meet $\Sigma_4$ at $A_4'$ then $A_4A_4'$ passes through $P$. Similarly if circle $CAD_4$ meets $\Sigma_4$ at $B_4'$ then $B_4B_4'$ passes through $P$ etc.

Altogether there must be twelve such lines passing through $P$, three for each of the circles $\Sigma_k$, $k$ = 1, 2, 3, 4. However because of the symmetry of the configuration only one such case needs to be established.

**Theorem 7**
The centres $Q_k$ of the circles $\Sigma_k$, $k$ = 1, 2, 3, 4 are collinear. Their radii are in proportion to their distances from $P$.

In Fig. 3 we illustrate all the above properties.

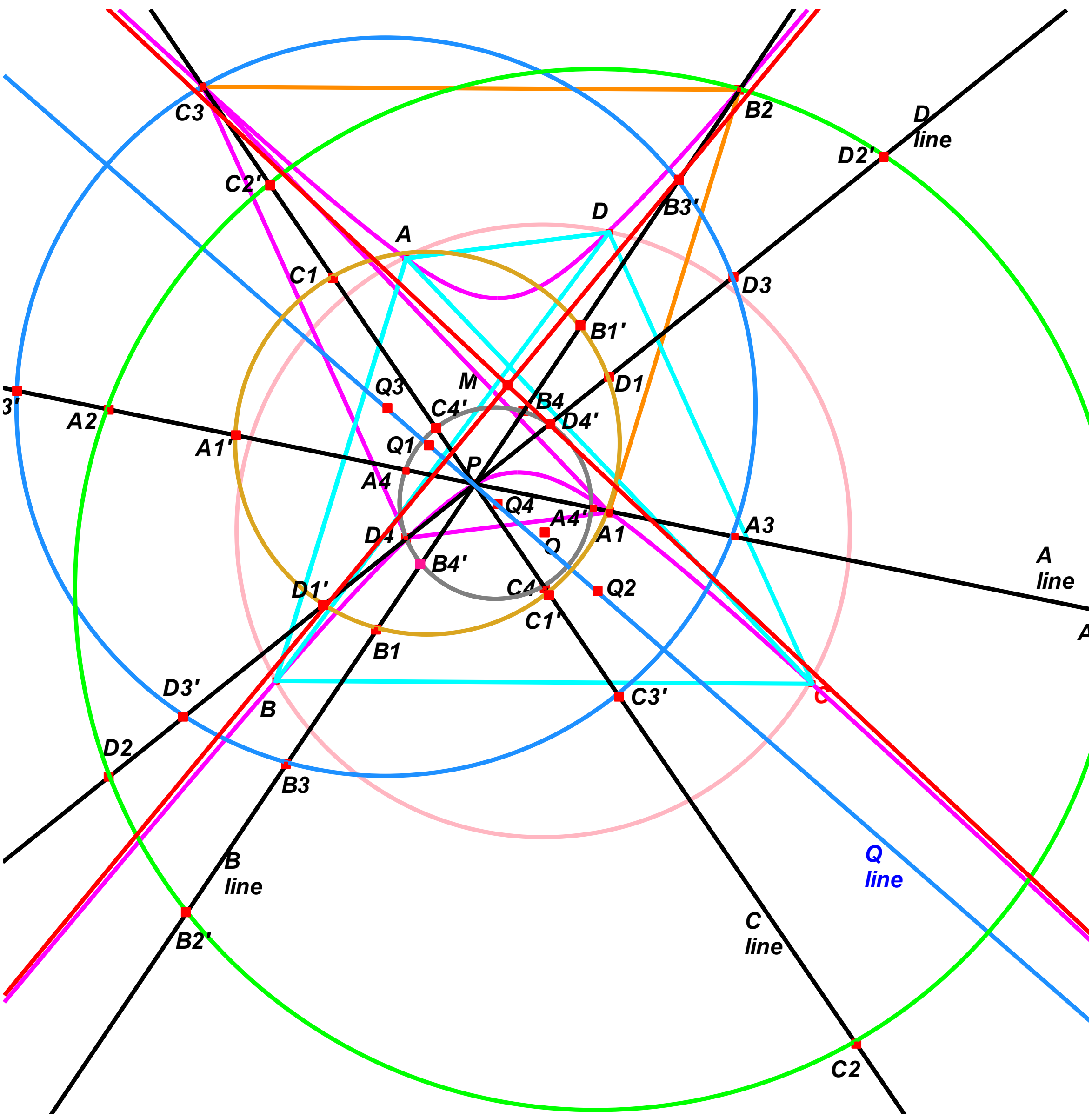

**Fig. 3**
**The A, B, C, D, Q lines**

## 6. Analysis and proofs

Consider the circle $\Gamma$ with equation

$$4m^2(x^2 + y^2) + m(m^2 + 1)(bdc + acd + abd + abc - a - b - c - d)x - (m^2 + 1)(bdc + acd + abd + abc + a + b + c + d)y + (m^2 + 1)(ab + ac + ad + bc + bd + cd) - 2(m^2 - 1) = 0. \quad (6.1)$$

It may be verified that it contains the points $(x, y)$ where

$$x = \{1/(2m)\}(t - 1/t),\ y = \{(1/2)\}(t + 1/t), \quad (6.2)$$

for $t = a, b, c, d$, provided $abcd = 1$ These define parametrically the hyperbola with equation $y^2 - m^2x^2 = 1$ and the points with parameters $a$, $b$, $c$, $d$ may be taken to define the points $A$, $B$, $C$, $D$, the cyclic quadrilateral inscribed in $\Gamma$. The centre $M$ of the hyperbola is the origin and the image of $ABCD$ under the 180° rotation about $M$ defines the congruent quadrilateral $A_1B_2C_3D_4$ and moreover the points $A_1$, $B_2$, $C_3$, $D_4$ lie on the hyperbola and have parameters $-a$, $-b$, $-c$, $-d$ respectively. The condition that four points on the hyperbola are concyclic is that the product of their parameters is 1. It follows that the sets of points $(ABC_3D_4)$, $(AB_2CD_4)$, $(AB_2C_3D)$, $(A_1BCD_4)$, $(A_1BC_3D)$, $(A_1B_2CD)$ are concyclic and that their equations may be obtained by altering signs in the equation of $\Gamma$ as appropriate. The equation of the chord of the hyperbola joining points with parameters $s$ and $t$ is

$$m(1 - st)x + (1 + st)y = s + t. \quad (6.3)$$

The point $P$ is now chosen on the hyperbola with parameter $p$. We can now put $(s, t)$ successively equal to $(a, -d)$, $(-a, p)$ and get the equations of the lines $D_4A$ and $A_1P$. Their intersection is by definition the point $A_4$, which has co-ordinates $(x, y)$, where

$$x = \{(a^2 + 1)(d + p) - 2a(dp + 1)\}/\{2am(p - d)\}, \quad (6.4)$$
$$y = \{(a^2 - 1)(d + p) + 2a(1 - dp)\}/\{2a(p - d)\}. \quad (6.5)$$

The co-ordinates of $B_4$, $C_4$ may now be obtained by using parameters $b$, $c$ rather than $a$. Points $B_1$, $C_1$, $D_1$ follow by using parameter $a$ instead of $d$ and $b$, $c$, $d$ respectively instead of $a$. Points $A_2$, $C_2$, $D_2$ follow by using parameter $b$ instead of $d$ and $a$, $c$, $d$ respectively instead of $a$. Points $A_3$, $B_3$, $D_3$ follow by using parameter $c$ instead of $d$ and $a$, $b$, $d$ respectively instead of $a$. The equation of the line $A_1P$ is worth recording as it is what we have termed the $A$-line. It has equation

$$m(1 + ap)x + y(1 - ap) = p - a. \quad (6.6)$$

**Proof of Theorem 5**

Take the co-ordinates of $A_4$, $B_4$, $C_4$, $D_4$ and construct a 4 x 4 matrix consisting of rows with entries $(x^2 + y^2, x, y, 1)$ for each of the four points. Then take its determinant and factorize and *DERIVE* provides the answer

$$\frac{(a-b)(a-c)(a-d)(b-c)(b-d)(c-d)(d+p)^4(1+m^2)(abcd-1)}{8a^2b^2c^2d^2(d-p)^4}. \qquad (6.7)$$

Since $abcd = 1$ it follows that the four points are concyclic on a circle we denote by $\Sigma_4$. A similar proof establishes the existence of circles $\Sigma_1$, $\Sigma_2$, $\Sigma_3$.

We now determine the equation of circle $\Sigma_4$. This is done using *Derive* by the same method as in the proof of Theorem 5, but by using current co-ordinates instead of those of $D_4$. The result is

$$4abcm^2(1-abcp)(x^2+y^2)+(ab+1)(ac+1)(bc+1)(abcp+1)m^3x+ \big(abcp(abc(abc+a+b+c)+bc+ca+ab-7)+(abc(a+b+c-7abc)+\ bc+ca+ ab+1)\big)mx+\{p\left(a^2b^2c^2\left(a+b+c-abc\right)-\ abc(bc+ca+ab+7)\right)+\ \left(\left(abc(a+b+ c+7abc)-(bc+ca+ab)+\ 1\right)\right)\Big\}m^2y+(1-ab)(1-bc)(1-ca)(1+abcp)y+ m^2\left(a^2b^2c^2\left(2abc-2p-p(bc+ca+ab)\right)+\ abc(a+b+c)p-abc(bc+ca+ab-2)+ a+b+c-2p\right)+\left(a^2b^2c^2\left(2abc+2p-p(bc+ca+ab)\right)+\ abc(a+b+c)p- abc(bc+ca+ab+2)+\ a+b+c-2p\right)=0 \qquad (6.8)$$

The line $AP$ meets $\Sigma_4$ at the point with co-ordinates $(x, y)$ where

$x = -(1/(2bcm(a^2p^2(m^2+1)+2ap(m^2-1)+m^2+1)))$ x $(a^2bcp(m^2+1)(b(c-p)-cp+1)- a(b^2c(c-p)(m^2+1)-b(c^2p(m^2+1)-c(3m^2(p^2-1)-p^2+1)-p(m^2+1))+p(c-p)(m^2+1))+ b(m^2+1)(1-cp)+c(m^2+1)-m^2p-p)$, (6.9)

$y = -(1/(2bc(a^2p^2(m^2+1)+2ap(m^2-1)+m^2+1)))$ x $(a^2bcp(m^2+1)(b(c-p)-cp-1)+a(b^2c(c-p)(m^2+1)-b(c^2p(m^2+1)+c(m^2(p^2+1)-3p^2-3)+p(m^2+1))-p(c-p)(m^2+1))-b(m^2+1)(cp+1)-c(m^2+1)+m^2p+p)$. (6.10)

When these co-ordinates are substituted into the equation of circle $BD_4A_1C$ they are found to satisfy it, and hence they are the co-ordinates of the point $A_4'$ lying on the $A$-line and the two circles. Similar analysis provides all twelve pairs of points on circles $\Sigma_1$, $\Sigma_2$, $\Sigma_3$, $\Sigma_4$ that are in involution through the involution point $P$. This completes the proof of Theorem 6.

The equations of the circles $\Sigma_k$, $k = 1, 2, 3, 4$ having been obtained it is now possible to find the co-ordinates of their centres. The co-ordinates of $Q_1$, $Q_4$, $P$ may then be used to show these points are collinear. This is done by forming the determinant whose three rows are $(x, y, 1)$ where $(x, y)$ are successively the co-ordinates of the three points. The value of this determinant turns out to be

$$(1/(32ab^2c^2d(abcp - 1)(bcdp - 1))) \text{ x } ((a - d)(m^2 + 1)(abcd - 1)(bcdp + 1)(b^3c^2p(m^2 + 1) + b^2c^2(cp(m^2 + 1) - 2(m^2(p^2 - 1) + p^2 + 1)) - b(2c(m^2(p^2 - 1) - p^2 - 1) + p(m^2 + 1)) - cp(m^2 + 1))). \tag{6.11}$$

This vanishes, on account of the factor $(abcd - 1)$ and so the three points are indeed collinear. Similarly any other two circle centres are collinear with $P$, so we have identified the existence of a Q-line containing all the circle centres and the point P. This establishes the first part of Theorem 7.

The second part of Theorem 7 is not difficult to check. The co-ordinates of all points are now known and it soon follows that triangles such as $A_4B_4C_4$ and $A_1B_1C_1$ are similar and that the enlargement factor is the same as $PQ_1/PQ_4$ (distances along the $Q$-line are signed, depending on which side of $P$ the centres are, and this corresponds to whether the triangles are directly homothetic through $P$ or whether a 180° twist is involved as well). See Fig. 3 again.

## 7. Results for the second generalization

Choose any point $T$, not on the sides of triangle $ABC$ or their extensions, which will act as a pseudo-orthocentre. In the original Hagge configuration $H$ acts as a centre of perspective and as an orthologic centre of triangles $ABC$ and $XYZ$. The generalization that is appropriate is to ensure that $T$ is an orthologic centre. An arbitrary circle is now drawn through $T$, which will serve as the generalized Hagge circle. *CABRI* indicates that there are generally two positions of $T$ on the circle for which it will also act as a centre of perspective, but the configuration shown in Fig. 4 is not one of these cases, as that is not necessary for the generalization.

To force $T$ to be an orthologic centre drop the perpendiculars from $T$ onto $BC$, $CA$, $AB$ and suppose these lines meet the circle through $T$ at points $X$, $Y$, $Z$ respectively. This forces $XYZ$ to have an orthologic centre with respect to $ABC$. Also the perpendiculars mean that, as may be proved by simple angle chasing, that triangle $XYZ$ is similar to $ABC$ and ordering of labels show that it is indirectly similar to $ABC$. It follows that the orthologic centre of $ABC$ with respect to $XYZ$ exists and is the point $J$ in the figure. $J$ and $T$ will be corresponding points in the indirect similarity that necessarily arises as a result of the orthologic property, so $J$ will lie on the circumcircle of $ABC$.

Standard methods may be used to find the centre *P* of inverse similarity and the axes of reflection (lines of inverse similarity). Only one is shown in Fig. 3.4.

Having located *P* draw *AP*, *BP*, *CP* to meet the circumcircle at *D*, *E*, *F* respectively. Find the images *U*, *V*, *W* under the indirect similarity. These are bound to lie on the circle through *T*. Also triangles *DEF* and *UVW* are bound to be similar, as they are related by the indirect similarity.

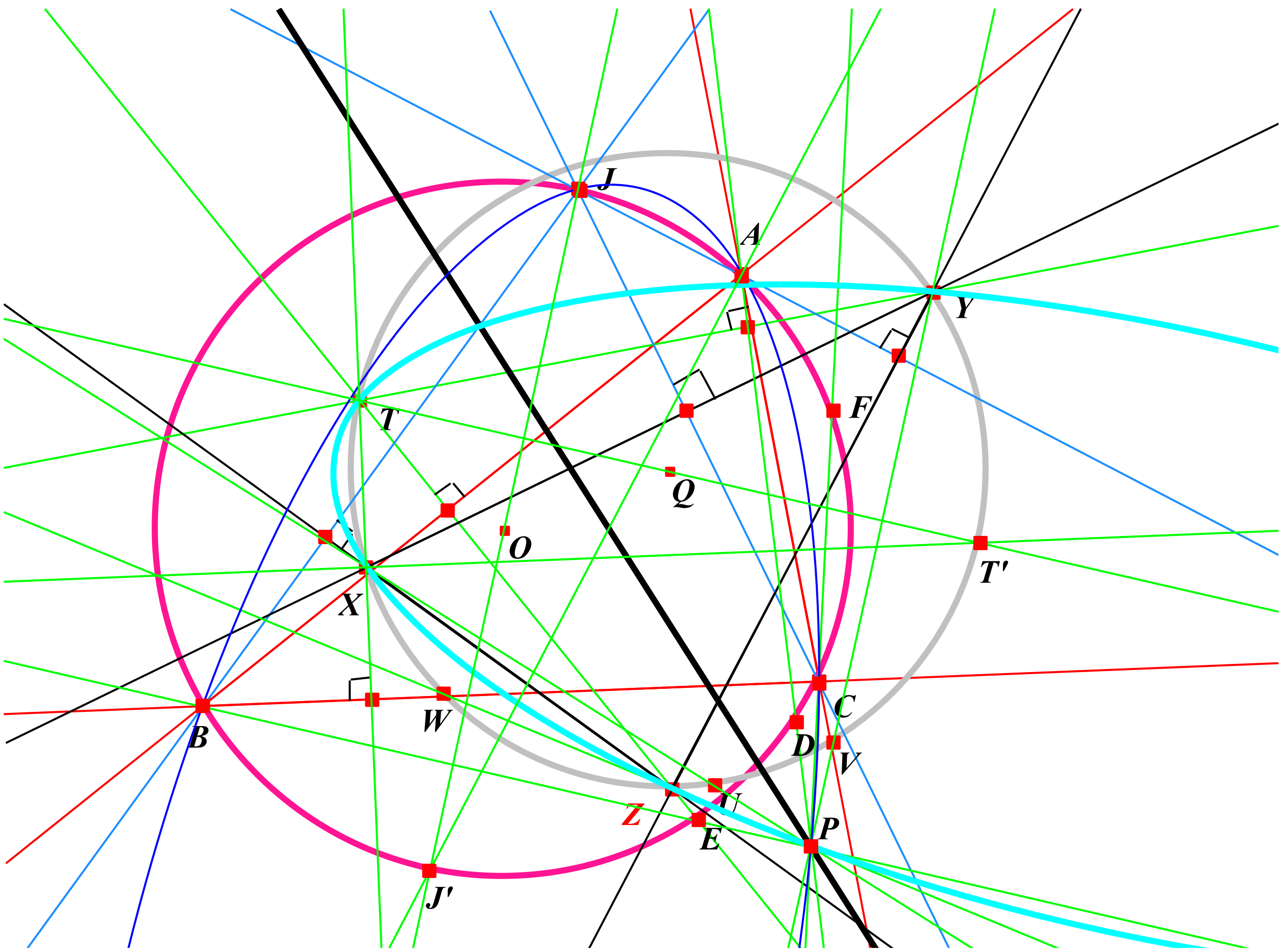

**Fig.4**
**Showing an indirect similarity and the orthologic centres**

The following theorem now holds:

**Theorem 8**

*XPU*, *YPV*, *ZPW* are straight lines. In other words triangles *XYZ*, *UVW* are in perspective with centre *P*, or since they all lie on a circle *X*, *U* and *Y*, *V* and *Z*, *W* are pairs in an involution on the circle through *T* by projection through *P*.

Despite there not being a perspective, some of Speckman's results still hold or generalize. In particular the conics *ABCPJ* and *XYZPT* will be images of each other in the indirect similarity. These conics are not necessarily hyperbolae; and even then *CABRI* confirms that their asymptotes are parallel only when a perspective exists between triangles *ABC* and *XYZ*.
However, triangles *ABC* and *XYZ* are paralogic, with paralogic centres at the opposite ends of the diameter to *T* and *J* respectively. These are labeled *T'* and *J'* in the figure.

**Theorem 9**

Let *VW* meet *XT* at *L*, with *M*, *N* similarly defined, then *LMN* is a straight line.

**Theorem 10**

The midpoint conic (midpoints of *AX*, *BY*, *CZ*, *DU*, *EV*, *FW*) exists and is illustrated in Fig. 5 along with the axis *LMN*.

Theorem 10 is, of course, a general theorem concerning six points on a pair of conics connected by an indirect similarity and needs no separate proof. Theorem 9 is unproved, but *CABRI* indicated. See Fig. 5.

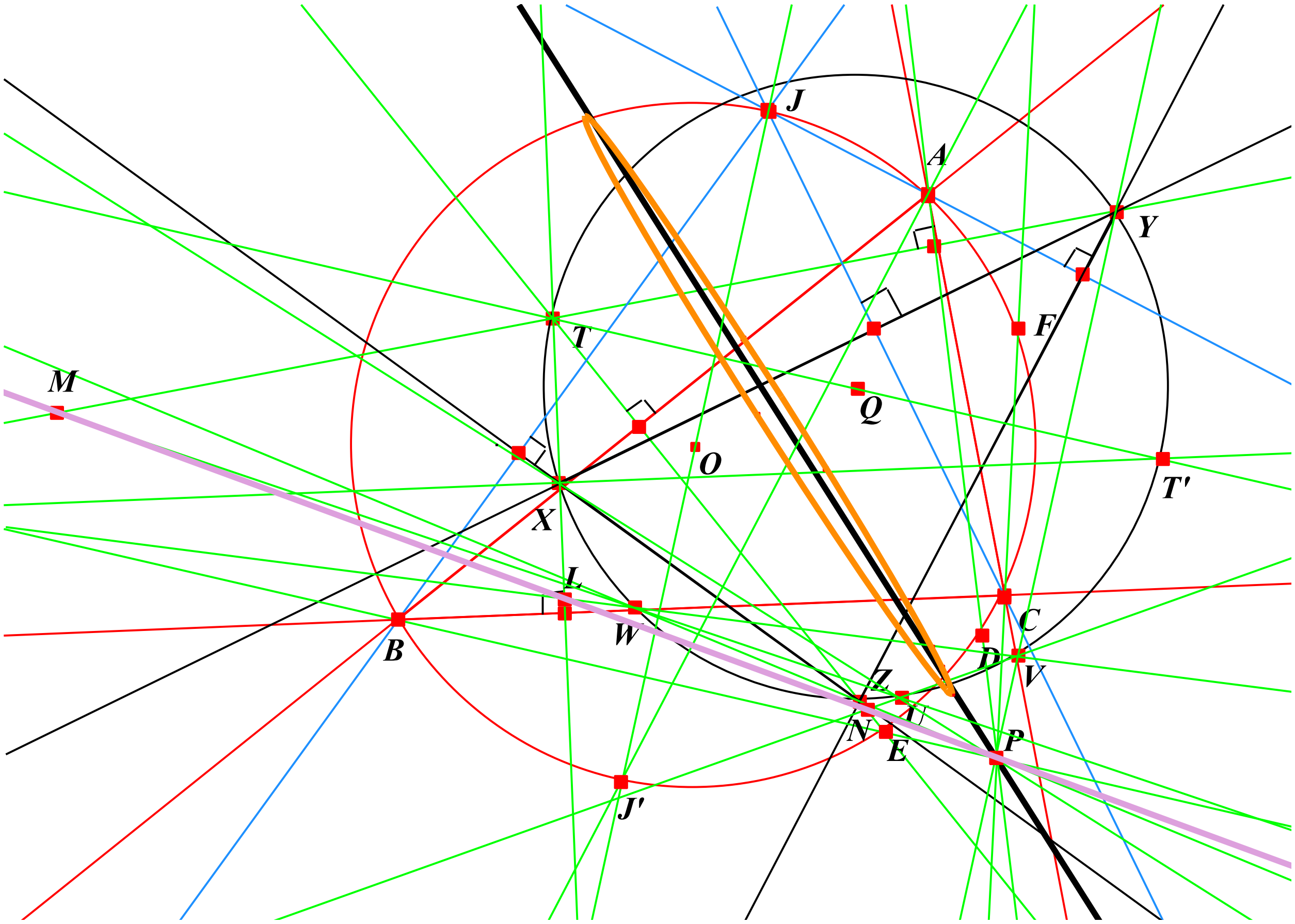

**Fig. 5**
**Cabri indicates the midpoint conic**

**Proof of Theorem 8**

This is trivial since it is a direct consequence of the indirect similarity. The points *A*, *P* and *D* are mapped by the indirect similarity on to points *X*, *P* and *U*. Since *APD* is a straight line, it follows that *XPU* is a straight line, *P* being the only invariant point and lines being mapped into lines. Similarly *YPV* and *ZPW* are straight lines.

*References*

1. A.M. Peiser, The Hagge circle of a triangle, *Amer. Math. Monthly*, 49 (1942) 524-527.
2. I. Boreico, *Result attributed to him* (Private communication).

Flat 4
Terrill Court
12-14 Apsley Road,
BRISTOL BS8 2SP.